\documentclass[11pt]{article}
\usepackage{amsfonts,graphics}

\usepackage{cite}
\usepackage{latexsym}
\usepackage{amsmath}
\usepackage{amssymb}
\usepackage[dvips]{epsfig}
\usepackage{amscd}

\newtheorem{theorem}{Theorem}[section]
\newtheorem{remark}{Remark}[section]

\newtheorem{lemma}{Lemma}[section]

\newcommand{\pa}{\partial}

\newcommand{\bi}{\bibitem}

\newcommand{\bt}{\begin{theorem}}
\newcommand{\bl}{\begin{lemma}}
\newcommand{\el}{\end{lemma}}
\newcommand{\et}{\end{theorem}}

\newcommand{\de}{\delta}
\newcommand{\ve}{\varepsilon}
\newcommand{\la}{\label}
\newcommand{\ka}{\kappa}

\newcommand{\ol}{\overline}

\newcommand{\bn}{\begin{eqnarray}}
\newcommand{\en}{\end{eqnarray}}
\newcommand{\bess}{\begin{eqnarray*}}
\newcommand{\eess}{\end{eqnarray*}}

\newcommand{\be}{\begin{equation}}
\newcommand{\ee}{\end{equation}}

\newcommand{\ba}{\begin{aligned}}
\newcommand{\ea}{\end{aligned}}

\def\norm[#1]#2{\|#2\|_{#1}}

\def\O{\Omega}
\def\rr{\mathbb{R}}

\makeatletter      
\@addtoreset{equation}{section}
\makeatother       
\date{}

\title{Some Uniform Estimates and Large-Time Behavior of Solutions to  One-Dimensional
Compressible Navier-Stokes
System in Unbounded Domains with Large Data\thanks{Partially supported by the
National Center for Mathematics and Interdisciplinary Sciences, CAS, and NNSFC
11371348, 11226163, and 11301422.}}
\date{}
\author{   Jing L{\small I}$^{a},$ Zhilei L{\small IANG}$^b$\thanks{
  Email addresses: ajingli@gmail.com (J. Li),   zhilei0592@gmail.com (Z. Liang)}
 \\[3mm] {\normalsize
  a. Institute of Applied Mathematics, AMSS,}\\
  {\normalsize \&   Hua Loo-Keng Key Laboratory of Mathematics,} \\
 {\normalsize Chinese Academy of Sciences,} \\
  {\normalsize Beijing 100190,  P. R. China } \\
 {\normalsize  b. School of Economic Mathematics, }\\ {\normalsize  Southwestern University of Finance  and
  Economics,} \\{\normalsize  Chengdu 611130, P.  R.
China}}

\begin{document}
 \maketitle



\begin{abstract}This paper is concerned with the large-time behavior of solutions
to the initial and initial boundary value problems with  large initial data for the compressible Navier-Stokes system  describing the  one-dimensional motion of a
viscous heat-conducting perfect polytropic  gas in unbounded domains.
The temperature is  proved to be  bounded  from below and above independently of both time and space. Moreover, it is shown that the global solution is    asymptotically stable as time tends to infinity.
Note that the  initial data can be arbitrarily large.
This result is proved by using elementary energy methods. \end{abstract}


\section{Introduction}

The compressible Navier-Stokes system  describing the  one-dimensional motion of a
viscous heat-conducting perfect polytropic  gas
can be written in the Lagrange variables in the following form (see \cite{ba,se})
\be \la{1.1}
v_t=u_{x},
\ee
\be \la{1.2}
u_{t}+P_{x}=\mu\left(\frac{u_{x}}{v}\right)_{x},\ee
\be  \la{1.'3} \left(e+\frac{u^2}{2}\right)_{t}+ (P
u)_{x}=\left(\kappa\frac{\theta_{x}}{v}+\mu\frac{uu_x}{v}\right)_{x}
 ,
\ee \be  \la{1.3'}   P =R \theta/{v},\quad e=c_v\theta +\mbox{ const.},
\ee where   $t>0$ is time,  $x\in\O\subset\rr=(-\infty,\infty)$ denotes the
Lagrange mass coordinate,   the unknown functions $v>0,u,$ $\theta>0,e>0,$ and $P$ are,  respectively, the specific volume of the gas,  fluid velocity,   internal energy,  absolute temperature, and  pressure,
$\mu $ and $\ka $ are the viscosity
and heat conductivity coefficients, $R>0$ is the gas constant, and $c_v $ is heat
capacity at constant volume. We assume
that $\mu,\ka,$ and $c_v$ are positive constants.

The system \eqref{1.1}-\eqref{1.3'} is supplemented
with the initial condition 
\be \la{1.4} (v(x,0),u(x,0),\theta(x,0))=(v_0(x),u_0(x),\theta_0(x)),\quad x\in
\O,\ee
and   three types of far-field and  boundary conditions:

1) Cauchy problem
  \be \la{1.5}  \O=\rr,\, \lim\limits_{|x|\rightarrow
  \infty}(v(x,t),u(x,t),\theta(x,t))=(1,0,1),\quad t>0;\ee

   2)   boundary and far-field conditions for $ \O=(0,\infty),$
\be \la{1.6}  u(0,t)=0,\, \theta_x(0,t)=0,\,\lim\limits_{x\rightarrow
\infty}(v(x,t),u(x,t),\theta(x,t))=(1,0,1),\quad t>0; \ee

 3)   boundary and far-field conditions for $ \O=(0,\infty),$
\be \la{1.7} u(0,t)=0,\, \theta (0,t)=1,\, \lim\limits_{x\rightarrow
\infty}(v(x,t),u(x,t),\theta(x,t))=(1,0,1), \quad t>0. \ee

 Kanel \cite{kan}    considered the Cauchy problem of the model system of equations
 \eqref{1.1} \eqref{1.2} with   $P=Rv^{-\gamma} $ and obtained   both the existence
 and the large-time asymptotic
behavior  of the  global  solutions    for  large initial data. For system
\eqref{1.1}-\eqref{1.3'},  Kazhikhov and Shelukhin
\cite{ks} first obtained    the global existence
     of solutions  in bounded domains  for large initial data. From then on,
significant progress has been made on the mathematical aspect of the initial and
initial
boundary value problems. For initial boundary value
problems in bounded
domains the existence and uniqueness of global (generalized) solutions and the
regularity
have been known. Moreover, the global solution is asymptotically stable as time
tends
to infinity; see \cite{akm,az1,az2,na1,na2,na3,ni,qy} among others.
For the Cauchy problem  \eqref{1.1}-\eqref{1.5} and the initial boundary value
problems
  \eqref{1.1}-\eqref{1.4}   \eqref{1.6}   and \eqref{1.1}-\eqref{1.4}   \eqref{1.7}
  (in unbounded domains), Kazhikhov \cite{ka1} (also cf. \cite{akm,ji4}) proved
  that
\begin{lemma}\la{q1k} Assume that the initial data $(v_0,u_0,\theta_0)$ satisfy
\be\la{c5} v_0-1,u_0,\theta_0-1\in H^1(\O),\,\, \inf_{x\in \O}v_0(x)>0,\, \,
\inf_{x\in \O}\theta_0(x)>0 ,\ee
and   are compatible with \eqref{1.6},  \eqref{1.7}.
Then there exists a unique global
(large)  generalized solution
  $(v ,u ,\theta )$   with positive  $v(x, t)$ and  $\theta(x, t)$   to
  \eqref{1.1}-\eqref{1.5}, or \eqref{1.1}-\eqref{1.4}  \eqref{1.6}, or
  \eqref{1.1}-\eqref{1.4}  \eqref{1.7}  satisfying that for any $T>0,$\be
  \la{3k}\begin{cases}   v-1,\,u,\,\theta-1 \in L^\infty(0,T;H^1(\O)),\, \, v_t\in
  L^\infty(0,T;L^2(\O)), \\ u_t,\,\theta_t,\,v_{xt},\,u_{xx},\,\theta_{xx} \,\in
  L^2(0,T;L^2(\O)).\end{cases}\ee
\end{lemma}

The asymptotic behavior as  $t\rightarrow \infty $ of the solution has been
studied under some smallness conditions on the initial data; see
\cite{ho,ji3,lz,ok,ka,knq,qy} and the references therein.
However, there are few results on the large-time behavior of the solution in the
case
of large data.  Jiang \cite{ji1,ji2} first obtained some interesting results on the
large-time behavior of solutions for large initial data by proving that
the specific volume is pointwise bounded from below and above independently of both time and space, and that   for all $t\ge 0$ the
temperature is bounded from below and above locally in $x.$ In particular, Jiang \cite{ji1,ji2} showed that\begin{lemma}[\cite{ji1,ji2}] \la{c4} Under the conditions of Lemma \ref{q1k}, let $(v ,u ,\theta
)$ be a generalized solution to \eqref{1.1}-\eqref{1.5}, or \eqref{1.1}-\eqref{1.4}
\eqref{1.6}, or \eqref{1.1}-\eqref{1.4}  \eqref{1.7} satisfying \eqref{3k} for any $T>0.$  Then there exists a positive constant $C_1$ depending only on
$\mu,\ka,R,c_v,$ $\|(v_0-1,u_0,\theta_0-1)\|_{H^1(\O)},
\inf\limits_{x\in \O}v_0(x),$ and $ \inf\limits_{x\in \O}\theta_0(x) $  such
that 
\be \la{c}C_1^{-1}\le v(x,t)\le C_1,\,\,\mbox{ for all } (x,t)\in \ol
\O\times [0,\infty).\ee\end{lemma}
From then on, for large initial data, whether the temperature  is pointwise bounded
from below and above independently of both time and space or not remains completely
open. This is an interesting problem partially because it is the key to study the large-time
dynamical behavior of the global generalized
solutions to \eqref{1.1}-\eqref{1.5},
  \eqref{1.1}-\eqref{1.4}   \eqref{1.6},  and \eqref{1.1}-\eqref{1.4}
  \eqref{1.7}. In   this paper, we will  give a positive answer and further prove that  the global solution is asymptotically stable as time
tends to infinity   for large initial data. Our main result is as follows:

\begin{theorem}\la{1k} Under the conditions of Lemma \ref{q1k}, let $(v ,u ,\theta
)$ be the (unique) generalized solution to \eqref{1.1}-\eqref{1.5}, or \eqref{1.1}-\eqref{1.4}
\eqref{1.6}, or \eqref{1.1}-\eqref{1.4}  \eqref{1.7} satisfying  \eqref{3k}  for any $T>0.$
  Then there exists a positive constant  $C_0$ depending only on
$\mu,\ka,R,c_v,$ $\|(v_0-1,u_0,\theta_0-1)\|_{H^1(\O)},
\inf\limits_{x\in \O}v_0(x),$ and $ \inf\limits_{x\in \O}\theta_0(x) $   such
that
\be \la{c2} C_0^{-1}\le \theta(x,t)\le C_0,\,\,\mbox{ for all } (x,t)\in \ol
\O\times [0,\infty) ,\ee
\be \sup\limits_{0\le t<\infty} \|(v-1,u,\theta-1)\|_{H^1(\O)} +\int_0^\infty
\left(\|v_x\|_{L^2(\O)}^2+\| (u_x,\theta_x)\|_{H^1(\O)}^2\right)dt\le C_0.\ee
Moreover, the following large-time behavior holds
\be\la{c3} \lim\limits_{t\rightarrow
\infty}\left(\|(v-1,u,\theta-1)(t)\|_{L^p(\O)}+
\|(v_x,u_x,\theta_x)(t)\|_{L^2(\O)}\right)=0,\ee for  any  $p\in (2,\infty].$

\end{theorem}

 \begin{remark} In Theorem \ref{1k}, we only assume that the initial data satisfy the  conditions which are needed for  the global
 existence of generalized solutions (see Lemma \ref{q1k}). Therefore, our results
 greatly improve the previous ones   due to \cite{ho,lz,ok,ka,knq,qy,ji3} where  some additional  smallness    conditions on    the initial data are
 needed.
\end{remark}

 \begin{remark} 
  For large initial data, Theorem
 \ref{1k}  shows that
the temperature is  bounded  from below and above independently of both time and
space and that the global solution   converges  to the constant steady
state  uniformly with respect to the spatial variable  as time goes to infinity.
Therefore, our results improve    those due to Jiang  \cite{ji1,ji2} where he proved  that  the temperature is  uniformly (in time)  bounded  from below and above   locally in $x$ and that global solutions are convergent locally in space  as time  goes to infinity.
\end{remark}

We now make some comments on the analysis of this paper. The key step to study the large-time behavior  of  the global
 generalized solutions is to get the   $L^2$-norm   (in both space and time)  bound of $\theta_x$ (see   \eqref{df8}). In fact,   \eqref{df8} has also been obtained under some additional smallness conditions on the initial data; see
\cite{ji3,lz,ok,ka,knq,qy} and the references therein. However, in our case, since the initial data may be arbitrarily large,      to obtain \eqref{df8}, some new ideas are needed. The key observations are as follows:    The combination of the standard   energetic  estimate (see \eqref{2.12}) with \eqref{c} shows that for $  \O_2(t) \triangleq\{x\in \O
|\,\theta(x,t)>2\} ,$  $$\int_0^\infty \int_{\O  \setminus\O_2(t)} \theta_x^2 dxdt  $$ is bounded. Hence, it suffices  to estimate the integral $$A \triangleq\int_0^\infty \int_{ \O_2(t)} \theta_x^2dxdt .$$ In fact, to estimate $A ,$ we multiply    the equation for the temperature by $(\theta-2)_+$ (see \eqref{1.1-1}). Then, to control the most difficult term  appearing in \eqref{1.1-1}, motivated by
 \cite{hxw}, we multiply  the equation for the velocity by $2u(\theta-2)_{+}$  (see \eqref{1.2-1}).    After some careful analysis on the integration  by parts over $\O_2(t)$ (see \eqref{df3}) and multiplying the equation for the velocity by $u^3,$ we finally find that $A $ can be controlled by  (see \eqref{lia5}) $$\int_0^\infty \sup\limits_{x\in\O}(\theta-3/2)_+^2 (x,t)dt,$$
 which in fact is bounded by $C(\ve)+\ve A $ for any $\ve>0$ (see \eqref{lia22}). These are the key to the proof of \eqref{df8}, and once that is obtained, the
proof follows in the same way as in   \cite{ji3,lz,ok,ka, knq,qy}. The whole procedure will be carried out in the next section.

\section{Proof of Theorem \ref{1k}}

We begin with the standard   energetic  estimate, which is motivated by the second
law of thermodynamics and embodies the dissipative effects of viscosity and thermal
diffusion.
\begin{lemma} \la{2k} It holds that
\be  \label{2.12}
 \ba
& \sup_{0\le t<\infty}\int_{\O}\left(\frac{1}{2} u^{2}+R(v-\ln
v-1)+c_{v}(\theta-\ln \theta-1)\right)  \\&\quad+
 \mu\int_{0}^{\infty}\int_{\O}
 \frac{u_{x}^{2}}{v\theta}+\kappa\int_{0}^{\infty}\int_{\O} \frac{\theta_{x}^{2}}{v\theta^{2}}
\le C,\ea\ee
where (and in what follows)   $C $ and  $C_i (i=2,\cdots,4)$ denote
generic positive constants
 depending only on $\mu,\ka,R,c_v,\|(v_0-1,u_0,\theta_0-1)\|_{H^1(\O)},
 \inf\limits_{x\in \O}v_0(x),$ and $ \inf\limits_{x\in \O}\theta_0(x).$

\end{lemma}

{\it Proof.} Using \eqref{1.1}, \eqref{1.2}, and \eqref{1.3'}, we rewrite
\eqref{1.'3} as \be  \la{1.3} c_{v}\theta_{t}+ R\frac\theta v
u_{x}=\kappa\left(\frac{\theta_{x}}{v}\right)_{x}+\mu\frac{u_{x}^{2}}{v}.
\ee  Multiplying  (\ref{1.1}) by
 $R(1- {v}^{-1})$,  (\ref{1.2})  by $
u$,   (\ref{1.3}) by
$  1- {\theta}^{-1} $,  and adding them altogether, we obtain
\bess \ba&(u^2/2+R(v-\ln v-1)+c_{v}(\theta-\ln
\theta-1))_t+\mu\frac{u^2_x}{v\theta}
+\kappa\frac{\theta_x^2}{v\theta^2}\\&=\left(\frac{\mu u u_x}{v}-\frac{R u
\theta}{v}\right)_x +Ru_x+\kappa\left((1-\theta^{-1})
\frac{\theta_x}{v}\right)_x,\ea\eess which together with \eqref{1.5} or \eqref{1.6}
or \eqref{1.7} yields \eqref{2.12}. We finish the proof of Lemma \ref{2k}.


Next, we derive the following $L^2$-norm (in both space and time) bounds of $\theta u_x$ and $\theta_x,$   which are  essential in our analysis.
\begin{lemma} \label{l3.2} There exists some positive constant $C $   such that for
any $T>0,$
\be  \la{df8}\sup_{0\le t\le T}\int_{\O} \left[(\theta-1)^{2}+u^{4}\right]+
\int_{0}^T\int_{\O} \left[(1+\theta+u^{2})u_{x}^{2}+
       \theta_{x}^{2}\right]
  \leq C. \ee
 \end{lemma}

\emph{Proof.} The proof of Lemma \ref{l3.2} will be divided into three steps.

{\it Step 1.}
First,   for $t\ge 0$ and $a>1,$ denoting   \bess \Omega_{a}(t)\triangleq\{x\in \O
|\,\theta(x,t)>a\} ,\eess we derive from
 \eqref{2.12} that \be \la{nep1}\sup\limits_{0\le t
 <\infty}\int_{\Omega_a(t)}\theta \le C(a)\sup\limits_{0\le t
 <\infty}\int_{\O}(\theta-\ln \theta-1)\le C(a). \ee

Next, integrating  \eqref{1.3} multiplied  by $(\theta-2)_{+}\triangleq \max\{\theta-2,0\}$ over $\O\times(0,T)$
gives
    \begin{equation}\begin{split} \label{1.1-1}  &\frac{c_{v}}{2} \int_{\O}
    (\theta-2)_{+}^{2}+\kappa\int_{0}^{T}\int_{\Omega_2(t)}\frac{ \theta_x^{2}}{v}
    \\  & = \frac{c_{v}}{2}\int_{\O} (\theta_{0}-2)_{+}^{2}-R  \int_{0}^{T}\int_{\O}\frac{\theta}{v}u_{x}(\theta-2)_{+} +\mu\int_{0}^{T}\int_{\O}
    \frac{u_{x}^{2}}{v}(\theta-2)_{+}. \end{split}\end{equation}
 To estimate the last term on the right hand side of \eqref{1.1-1}, motivated by
 \cite{hxw},  we multiply   \eqref{1.2} by $2u(\theta-2)_{+}$ and integrate the
 resulting equality over $\O\times (0,T)$ to get
 \be \ba\begin{split} \label{1.2-1}&  \int_{\O}  u^2  (\theta-2)_{+}
+2\mu\int_{0}^{T}\int_{\O} \frac{u_{x}^{2}}{v}(\theta-2)_{+}\\&=\int_{\O}   u_0^2
(\theta_0-2)_{+}+2R \int_{0}^{T}\int_{\O} \frac{\theta}{v} u_{x}(\theta-2)_{+}
+ 2R \int_{0}^{T}\int_{\Omega_2(t)}\frac{\theta}{v} u \theta_x\\&\quad
 - 2\mu\int_{0}^{T}\int_{\Omega_2(t)}\frac{u_{x}}{v}u  \theta_x
 +\int_{0}^{T}\int_{\Omega_2(t)} u^{2} \theta_t. \end{split}\ea\ee
Adding  (\ref{1.2-1}) to (\ref{1.1-1}), we obtain after using \eqref{1.3} that
     \be \begin{split} \label{1.3-1}
 &\int_{\O}\left[\frac{c_{v}}{2} (\theta-2)_{+}^{2}+  u^{2}(\theta-2)_{+}\right]
 +\kappa\int_{0}^{T}\int_{\Omega_2(t)}\frac{\theta_x^{2}}{v}
 +\mu\int_{0}^{T}\int_{\O}\frac{u_{x}^{2}}{v}(\theta-2)_{+}
 \\
     &= \int_{\O}\left[\frac{c_{v}}{2} (\theta_0-2)_{+}^{2}+
     u_0^{2}(\theta_0-2)_{+}\right] +R\int_{0}^{T}\int_{\O}
     \frac{\theta}{v}u_{x}(\theta-2)_{+}\\
&\quad + 2R\int_{0}^{T}\int_{\Omega_2(t)} \frac{\theta}{v} u \theta_x
 - 2\mu\int_{0}^{T}\int_{\Omega_2(t)}\frac{u_{x}}{v}u  \theta_x \\
     &\quad  +\frac{1}{c_{v}} \int_{0}^{T}\int_{\Omega_{2}(t)}u^{2}
 \left(\mu\frac{u_{x}^{2}}{v}-R\frac{\theta}{v}u_{x}\right)+\frac{\ka}{c_{v}} \int_{0}^{T}\int_{\Omega_{2}(t)}u^{2}
    \left(\frac{\theta_{x}}{v}\right)_{x}\\
  &\triangleq \int_{\O}\left[\frac{c_{v}}{2} (\theta_0-2)_{+}^{2}+
  u_0^{2}(\theta_0-2)_{+}\right] +\sum_{i=1}^{5}I_{i}. \end{split}\ee

We estimate each $I_i(i=1,\cdots,5)$ as follows:

First, it follows from Cauchy's inequality  and \eqref{c}
that\be \ba\la{df1} |I_1|&=R \left|\int_{0}^{T}\int_{\O}
\frac{\theta}{v}u_{x}(\theta-2)_{+}\right|\\
 &\le
 \frac{\mu}{2}\int_{0}^{T}\int_{\O}\frac{u_{x}^{2}}{v}(\theta-2)_{+}+C\int_{0}^{T} \int_{\O}\theta^{2}(\theta-2)_{+}\\
 &\le  \frac{\mu}{2}\int_{0}^{T}\int_{\O}\frac{u_{x}^{2}}{v}(\theta-2)_{+} +C\int_{0}^{T}\int_{\O } \theta(\theta-3/2)^{2}_{+}\\
 &\le  \frac{\mu}{2}\int_{0}^{T}\int_{\O}\frac{u_{x}^{2}}{v}(\theta-2)_{+} +C\int_{0}^{T}\sup_{x\in
 \Omega }(\theta-3/2)^{2}_{+}(x,t),\ea\ee   where in the last inequality we have used  \eqref{nep1}.

Next, Cauchy's inequality  and \eqref{c}  yield that for any $\ve>0,$  \be \ba\la{df2}
 |I_{2}|+|I_3|&=2R\left|\int_{0}^{T}\int_{\Omega_2(t)}
 \frac{\theta}{v}u\theta_x\right|+2\mu
 \left|\int_{0}^{T}\int_{\Omega_2(t)}\frac{u_{x}}{v}u \theta_x\right|\\
 &\le \ve\int_{0}^{T}\int_{\O}\theta_x^2+C(\ve)\int_{0}^{T}
 \int_{\Omega_{2}(t)}u^{2}\theta^{2} +C(\ve)\int_{0}^{T}\int_{\O}u^{2}u_{x}^{2} \\
 &\le \ve\int_{0}^{T}\int_{\O}\theta_x^2
 +C(\ve)\int_{0}^{T}\sup_{x\in
 \Omega }(\theta-3/2)^{2}_{+}(x,t) +C(\ve)\int_{0}^{T}
 \int_{\O}u^{2}u_{x}^{2} , \ea\ee where in the last inequality we have used\be  \la{2.13}\int_{0}^{T}
 \int_{\Omega_{2}(t)}u^{2}\theta^{2} \le 16\int_{0}^{T}
 \int_{\Omega }u^{2}(\theta-3/2)^{2}_{+} \le C  \int_{0}^{T}\sup_{x\in
 \Omega }(\theta-3/2)^{2}_{+}(x,t), \ee due to \eqref{2.12}.

 Then,  it follows from Cauchy's inequality  and \eqref{2.13} that
\be\la{df}\ba
|I_4| &\le  C\int_{0}^{T}\int_{\O}u^{2}u_{x}^{2}+ C \int_{0}^{T}\sup_{x\in
 \Omega }(\theta-3/2)^{2}_{+}(x,t) .
\ea\ee

 Finally, for
\bess\ba   \varphi_{\eta}(\theta)\triangleq
                       \begin{cases}
                        1, & \theta-2>\eta, \\
                        (\theta-2)/\eta, & 0\le \theta-2\le \eta,\\
                         0, &  \theta-2\le 0,
                       \end{cases}
                     \ea\eess
Lebesgue's dominated convergence theorem shows that for any $\ve>0,$ \be \ba\la{df3} I_{5} &=
\frac{\ka}{c_{v}}\lim_{\eta\rightarrow0+}\int_{0}^{T} \int_{\O}
\varphi_{\eta} (\theta) u^{2}\left(\frac{\theta_{x}}{v}\right)_{x} \\
&= \frac{\ka}{c_{v}}\lim_{\eta\rightarrow0+} \int_{0}^{T}\int_{\O} \left(-2 \varphi_{\eta}(\theta ) uu_{x}\frac{\theta_{x}}{v}-
\varphi_{\eta}' (\theta ) u^{2}\frac{\theta_{x}^{2}}{v} \right)\\
&\le -\frac{2\ka}{c_{v}}  \int_{0}^{T}\int_{\O_2(t)}   uu_{x}\frac{\theta_{x}}{v} \\
&\le \ve\int_{0}^{T}\int_{\O}  \theta_x^{2}
+C(\ve)\int_{0}^{T}\int_{\O}u^{2}u_{x}^{2} ,
\ea\ee where in the third inequality we have used
 $\varphi_\eta'(\theta)\ge 0.$

Noticing that
\bess\ba  &\int_{0}^{T}\int_{\O}\left(u_{x}^{2}\theta+\theta_{x}^{2}\right)\\
&=\int_{0}^{T}\int_{\Omega_{3}(t)}\left(u_{x}^{2}\theta+\theta_{x}^{2}\right)+\int_{0}^{T}\int_{\O\setminus\Omega_{3}(t)}\left(u_{x}^{2}\theta+\theta_{x}^{2}\right)\\
&\le 3\int_{0}^{T}\int_{\Omega_{3}(t)}\left( u_{x}^{2}(\theta-2)_{+}+
{\theta_{x}^{2}} \right)
+C\int_{0}^{T}\int_{\O\setminus\Omega_{3}(t)}\left(\mu\frac{u_{x}^{2}}{\theta}
+\ka\frac{\theta_{x}^{2}}{\theta^{2}}\right)\\
&\le C\int_{0}^{T}\int_{\Omega_{2}(t)}\left(u_{x}^{2}(\theta-2)_{+}+\frac{
\theta_x^{2}}{v}\right)
+C,\ea\eess where in the last inequality we have used $\O_3(t)\subset\O_2(t),$ \eqref{c}, and \eqref{2.12},
we substitute \eqref{df1}, \eqref{df2}, \eqref{df}, and \eqref{df3} into \eqref{1.3-1} and choose $\ve$ suitably
small to obtain
    \be \begin{split} \label{df4}
 &\sup_{0\le t\le T}\int_{\O}
 (\theta-2)_{+}^{2}+\int_{0}^{T}\int_{\O}\left(u_{x}^{2}\theta+\theta_{x}^{2}\right)\\
&\le  C+ C  \int_{0}^{T}\sup_{x\in
 \Omega }(\theta-3/2)^{2}_{+}(x,t) +
C_2\int_{0}^{T}\int_{\O}u^{2}u_{x}^{2}. \end{split}\ee

{\it Step 2.} To estimate the last term on the right hand side of \eqref{df4}, we multiply
 \eqref{1.2} by $u^{3}$ and integrate the resulting equality over $\O\times (0,T)$
 to get
     \be \ba\la{lia}  &\frac{1}{4}\int_{\O}  u^{4}+3\mu
     \int_{0}^{T}\int_{\O} \frac{u^{2}u_{x}^{2}}{v} \\
&=\frac{1}{4}\int_{\O}  u^{4}_{0}+3R\int_{0}^{T}\int_{\O} \frac{ 1-v
}{v}u^{2}u_{x}+3R\int_{0}^{T}\int_{\O\setminus\O_2(t)} \frac{\theta-1}{v} u^{2}u_{x}\\&\quad+3R\int_{0}^{T}\int_{\O_2(t)} \frac{\theta-1}{v} u^{2}u_{x} \triangleq \frac{1}{4}\int_{\O}  u^{4}_{0}+\sum\limits_{i=1}^3 J_i.\ea\ee

  It follows from \eqref{2.12}
and \eqref{c} that for any $\alpha\in [2,3],$\be
\la{4y} \ba&\sup\limits_{0\le t\le T}\int_\O  (v-1)^2+\sup_{0\le t\le T}\int_{\O\setminus
\Omega_{\alpha}(t)}(\theta-1)^{2} \\&\le  C\sup_{0\le t\le T}\int_{\O}(v-\ln v-1)+ C \sup_{0\le t\le T}\int_{\O}(\theta-\ln \theta-1) \le C ,\ea\ee
 which together with Holder's inequality yields that
  \be\la{cc}\ba |J_1|+|J_2|&\le C\int_0^T\sup\limits_{x\in\O}  u^2(x,t) \|u_x \|_{L^2(\O)}\left(\int_\O  (v-1)^2+ \int_{\O\setminus
\Omega_{2}(t)}(\theta-1)^{2}\right)^{1/2}\\&\le C\int_0^T \int_\O u_x^2  , \ea\ee where in the second inequality we have used   \eqref{2.12} and the following simple fact that for any $w\in H^1(\O),$
  \be  \la{pj}\ba \sup\limits_{x\in\O} w^2(x ) &=\sup\limits_{x\in\O}\left(-2 \int_x^\infty w(y )w_x(y )dy\right)\\&\le 2\|w \|_{L^2(\O)}\|w_x \|_{L^2(\O)}. \ea\ee
The combination of Cauchy's inequality with \eqref{2.13} leads to
  \be \ba\la{rf1}
|J_3|&\le  \mu\int_{0}^{T}\int_{\Omega_{2}(t)} \frac{u^{2}u_{x}^{2}}{v}
+C\int_{0}^{T}\int_{\Omega_{2}(t)}  \theta^{2}u^{2}  \\
&\le\mu \int_{0}^{T}\int_{\O}  \frac{u^{2}u_{x}^{2}}{v}+ C\int_{0}^{T}
  \sup_{x\in
 \Omega }(\theta-3/2)^{2}_{+}(x,t) .
\ea\ee
Putting   \eqref{cc} and \eqref{rf1} into \eqref{lia}  gives
\be \ba\la{1.5-1}
 & \sup_{0\le t\le T}\int_{\O}  u^{4}+\int_{0}^{T}\int_{\O} u^{2}u_{x}^{2}\\
  &\le C+C  \int_{0}^{T}\int_{\O} u_{x}^{2}+C \int_{0}^{T}
 \sup_{x\in
 \Omega }(\theta-3/2)^{2}_{+}(x,t)   \\
  &\le C(\de) +C\delta \int_{0}^{T}\int_{\O} \theta u_{x}^{2}  +C\int_{0}^{T}
    \sup_{x\in  \Omega }(\theta-3/2)^{2}_{+}(x,t) ,
    \ea\ee where in the last inequality we have used the following simple fact that for any $\de>0,$\be\la{1pp} 2\int_{0}^{T}\int_{\O}u_x^2\le \de\int_{0}^{T}\int_{\O}\theta u_x^2  +  \de^{-1} \int_{0}^{T}\int_{\O}\theta^{-1} u_x^2 \le  \de \int_{0}^{T}\int_{\O}\theta u_x^2  +C (\de), \ee due to Cauchy's inequality, \eqref{2.12}, and \eqref{c}.

Adding \eqref{1.5-1}  multiplied  by $  C_2+1$ to   \eqref{df4}, then choosing
$\delta $ suitably  small, we have
\be \ba\label{lia5}
&\sup_{0\le t\le T}\int_{\O} \left[(\theta-2)_{+}^{2}+ u^{4}\right]+
\int_{0}^{T}\int_{\O}\left[(\theta+u^{2})u_{x}^{2}+
       \theta_{x}^{2}\right]\\
&\leq  C +C \int_{0}^{T}  \sup_{x\in
 \Omega }(\theta-3/2)^{2}_{+}(x,t) .
\ea\ee

{\it Step 3.} It remains to estimate the last   term  on the right hand side of
 \eqref{lia5}.
In fact, standard calculations yield   that for    any $\ve>0,$
\be \la{lia22} \ba
\int_0^T\sup\limits_{x\in\O}(\theta(x,t)-3/2)_{+}^{2}
&=\int_0^T\sup\limits_{x\in\O}\left(  \int^{\infty}_{x}\pa_x(\theta-3/2)_{+}  \right)^2\\
&\le \int_0^T \left(  \int_{\O_{3/2}(t)} |\theta_x|  \right)^2\\
&\le \int_0^T\left( \int_{\O_{3/2}(t)} \frac{ \theta_x^2}{\theta}\int_{\O_{3/2}(t)}  {\theta}\right)\\& \le C \int_0^T
\int_{\O}\frac{\theta_{x}^{2}}{\theta}\\ & \le C(\ve)\int_0^T
\int_{\O}\frac{\theta_{x}^{2}}{v\theta^2}+ \ve \int_0^T \int_{\O} {\theta_{x}^{2}}\\ & \le C(\ve) + \ve \int_0^T \int_{\O} {\theta_{x}^{2}}, \ea\ee where in the fourth and  last inequalities we have used \eqref{nep1} and  \eqref{2.12} respectively.
Putting \eqref{lia22}   into  \eqref{lia5} and choosing $\ve$ suitably small lead  to
  \bess   \sup_{0\le t\le T}
 \int_{\O} \left[(\theta-2)_{+}^{2}+
 u^{4}\right]+\int_{0}^{T}\int_{\O}\left[(\theta+ u^{2})u^{2}_{x}+
 \theta_{x}^{2}\right] \leq   C ,
 \eess   which combined with \eqref{4y} and \eqref{1pp}  immediately gives \eqref{df8}. The proof of Lemma \ref{l3.2} is
completed.

 We will  derive some necessary  uniform estimates on the spatial derivatives  of
the global generalized solution $(v,u,\theta)$ in the next lemma.
\begin{lemma} \la{5y}  There exists some positive constant $C $   such that for any
$T>0,$ \be\la{z4} \sup_{0\le t\le T}\int_{\O} \left(v_x^2+u_x^2+\theta_x^2\right)
+\int_0^T\int_{\O}\left(\theta v_x^2+u_{xx}^2+\theta_{xx}^2\right)\le
C.\ee\end{lemma}

{\it Proof.} First, integrating  \eqref{1.2}     multiplied by  $\frac{v_{x}}{v}$
over $\O ,$ we obtain after  using \eqref{1.1}  that
\bess\begin{split} \frac{\mu}{2}\frac{d}{dt}\int_{\O} \frac{v_{x}^2}{v^2}
&= R \int_{\O} \left(\frac{\theta}{v}\right)_{x}\frac{v_{x}}{v}+ \int_{\O} u_t
\frac{v_{x}}{v} \\&= R\int_{\O} \frac{\theta_{x}v_{x}}{v^{2}}-R \int_{\O}
\frac{\theta v^{2}_{x}}{v^{3}} + \frac{d}{dt}\int_{\O} u \frac{v_{x}}{v}-\int_{\O}
u_x\frac{ v_t}{v}\\&\le C\int_{\O}  \frac{\theta_{x}^2}{v\theta}- \frac{R}{2}\int_{\O}
\frac{\theta v^{2}_{x}}{v^{3}} + \frac{d}{dt}\int_{\O} u \frac{v_{x}}{v}+\int_{\O}
\frac{u_x^2}{v},\end{split}\eess
which together with Cauchy's inequality, \eqref{c},  \eqref{2.12},  \eqref{df8}, and \eqref{lia22}    gives
   \be \ba\la{5.5} \sup\limits_{0\le t\le T}\int_{\O}
   v_{x}^{2}+\int_{0}^{T}\int_{\O} \theta v_{x}^{2}
&\le C.
\ea\ee

Next, integrating \eqref{1.2} multiplied by $u_{xx}$ over $\O $ yields
\be \ba\la{z}&\frac{1}{2}\frac{d}{dt }\int_{\O}u_{x}^{2}+\mu
\int_{\O}\frac{u_{xx}^{2}}{v}=R\int_{\O}\left(\frac{\theta}{v}\right)_{x}u_{xx}+\mu
\int_{\O}\frac{u_{x}}{v^{2}}v_{x}u_{xx} .\ea\ee
It follows from  \eqref{df8},     \eqref{5.5}, and    \eqref{pj}   that
\be \ba\la{zz}& \int_0^T\left|R\int_{\O}\left(\frac{\theta}{v}\right)_{x}u_{xx}+\mu
\int_{\O}\frac{u_{x}}{v^{2}}v_{x}u_{xx}\right|  \\&\le
\frac{\mu}{4}\int_{0}^{T}\int_{\O}\frac{u_{xx}^{2}}{v}
+C\int_{0}^{T}\int_{\O}\left(\theta^{2}v_{x}^{2} +\theta^{2}_{x}+u_{x}^{2}v_{x}^{2}
\right) \\&\le C+ \frac{\mu}{4}\int_{0}^{T}\int_{\O}\frac{u_{xx}^{2}}{v}
+C\max_{\ol\O\times[0,T]}\theta\int_{0}^{T}\int_{\O}\theta v_{x}^{2}
+C\int_{0}^{T}\|u_{x}(\cdot,t)\|_{L^{\infty}(\O)}^{2}   \\&\le C+
\frac{\mu}{2}\int_{0}^{T}\int_{\O}\frac{u_{xx}^{2}}{v} +C\max_{\ol\O\times[0,T]}
\theta,\ea\ee
 which combined with \eqref{z} shows
 \be \ba\la{z1}& \sup_{0\le t\le T} \int_{\O}u_{x}^{2}+ \int_{0}^{T}\int_{\O}
 u_{xx}^{2} \le C+C\max_{\ol\O\times[0,T]} \theta.\ea\ee

 Next, integrating \eqref{1.3} multiplied by $\theta_{xx}$ over $\O $ leads to
\be  \la{nep7}\ba\frac{1}{2}\frac{d}{dt }\int_{\O}\theta_{x}^{2}
+\ka \int_{\O}\frac{\theta_{xx}^{2}}{v}
&= \ka \int_{\O}\frac{\theta_{x}v_{x}}{v^{2}}\theta_{xx} -\mu
\int_{\O}\frac{u_{x}^{2}}{v}\theta_{xx}+  \int_{\O} R\frac{\theta}{v}u_{x}
\theta_{xx} .\ea\ee
Cauchy's inequality and \eqref{pj} give
\be  \la{np8}\ba &\int_0^T\left|\ka
\int_{\O}\frac{\theta_{x}v_{x}}{v^{2}}\theta_{xx} -\mu
\int_{\O}\frac{u_{x}^{2}}{v}\theta_{xx}+  \int_{\O} R\frac{\theta}{v}u_{x}
\theta_{xx}\right| \\
&\le C\int_{0}^{T} \| \theta_{xx} \|_{L^2(\O)}  \|\theta_{x}\|_{L^\infty(\O)}
\|v_{x}\|_{L^{2}(\O)} \\&\quad+C\int_{0}^{T} \| \theta_{xx} \|_{L^2(\O)} \left(
  \|u_{x}\|_{L^\infty(\O)} \|u_{x}\|_{L^{2}(\O)}+  \|\theta\|_{L^\infty(\O)} \|
  u_{x}\|_{L^2(\O)} \right)  \\
&\le C\int_{0}^{T} \| \theta_{xx} \|_{L^2(\O)}
\|\theta_{xx}\|^{1/2}_{L^2(\O)}\|\theta_{x}\|^{1/2}_{L^2(\O)}\|v_{x}\|_{L^{2}(\O)}
 \\&\quad+C\int_{0}^{T}\| \theta_{xx} \|_{L^2(\O)}
 \|u_{x}\|_{H^1(\O)}\left(\|u_{x}\|_{L^{2}(\O)}+ \|\theta\|_{L^\infty(\O)} \right)
 \\
&\le \frac{\ka}{4}\int_{0}^{T}\int_{\O}\frac{\theta_{xx}^{2}}{v} + C
+C\max_{\ol\O\times[0,T]}\theta^3,\ea\ee where in the last inequality we have used
\eqref{5.5},  \eqref{df8}, and \eqref{z1}.
  Integrating \eqref{nep7} over $(0,T) ,$  we obtain after using \eqref{np8} that
 \be \ba\la{z3}& \sup\limits_{0\le t\le T}\int_{\O}\theta_{x}^{2}+
 \int_{0}^{T}\int_{\O} \theta_{xx}^{2} \le C+C\max_{\ol\O\times[0,T]}
 \theta^{3}.\ea\ee

Finally, it follows from \eqref{pj} and \eqref{df8} that for all $t\ge 0,$
 \be\la{z7}\ba  \|(\theta-1)(\cdot,t)\|_{C(\ol\O)}^{2}&\le C
 \|(\theta-1)(\cdot,t)\|_{L^{2}(\O)} \|\theta_{x}(\cdot,t)\|_{L^{2}(\O)}\\&\le C
 \|\theta_{x}(\cdot,t)\|_{L^{2}(\O)},\ea\ee
 which combined with \eqref{z3} yields
 \bess \max_{\ol\O\times[0,T]} (\theta-1)^{2}\le
 C+C\max_{\ol\O\times[0,T]}\theta^{3/2}.\eess
This implies that there exists a positive constant $C_3$ such that for any $(x,t)\in \ol\O\times [0,T],$
\be\la{np9}\theta(x,t)\le C_3, \ee
which together with \eqref{z1}, \eqref{z3}, and \eqref{5.5} gives  \eqref{z4} and
finishes the proof of Lemma \ref{5y}.

With  Lemma \ref{5y} at hand, we are now in a position to prove the following
  large-time behavior of   global generalized solutions which together with  Lemmas \ref{2k}-\ref{5y}   finishes the proof of Theorem \ref{1k}.
\begin{lemma}  \la{6y} It holds that \be \la{1y}\lim\limits_{t\rightarrow
\infty}\left(\|(v-1,u,\theta-1)(t)\|_{L^p(\O)}+
\|(v_x,u_x,\theta_x)(t)\|_{L^2(\O)}\right)=0,\ee  for  any   $p\in (2,\infty].$
Moreover, there exists a positive constant $C_4 $ such that  for all $(x,t)\in \ol \O\times
[0,\infty)$  \be \la{2y}C_4^{-1}\le \theta(x,t)\le C_4.\ee\end{lemma}

{\it Proof.} It follows from \eqref{df8},    \eqref{z},  \eqref{zz},
\eqref{nep7},  \eqref{np8},    \eqref{np9},  and \eqref{z4} that \bess \ba&
\int_0^\infty\left(\|u_x(\cdot,t)\|_{L^2(\O)}^2+
\left|\frac{d}{dt}\|u_x(\cdot,t)\|_{L^2(\O)}^2   \right|\right)dt\\&+
\int_0^\infty\left(\|\theta_x(\cdot,t)\|_{L^2(\O)}^2+
\left|\frac{d}{dt}\|\theta_x(\cdot,t)\|_{L^2(\O)}^2   \right|\right)dt \le C,\ea\eess
which directly gives \be \la{x1}\lim\limits_{t\rightarrow
\infty}\left(\|u_x(\cdot,t)\|_{L^2(\O)}+\|\theta_x(\cdot,t)\|_{L^2(\O)}\right)=0.\ee
This combined with \eqref{z7} shows \bess \lim\limits_{t\rightarrow \infty}\|\theta
(\cdot,t)-1\|_{C(\ol\O)}=0.\eess   Hence, there exists some $T_0>0$ such that for
all $(x,t)\in \ol \O\times [T_0,\infty)$ \be  \la{x5} 1/2\le \theta(x,t)\le 3/2,\ee
which, along with   \eqref{z4}, leads to \be  \la{x4} \int_{T_0}^\infty
\|v_x(\cdot,t)\|_{L^2(\O)}^2   \le C.\ee
This combined with \eqref{1.1} and \eqref{z4} yields
\bess   \ba\int_{T_0}^\infty \left|\frac{d}{dt}\|v_x(\cdot,t)\|_{L^2(\O)}^2
\right|   &=2\int_{T_0}^\infty \left|\int_\O u_{xx}v_x \right|  \\& \le
\int_{T_0}^\infty  \int_\O u_{xx}^2   +\int_{T_0}^\infty  \int_\O v_{x }^2 \le C
,\ea \eess
 which together with \eqref{x4} implies \be  \la{x8}\lim\limits_{t\rightarrow
 \infty}\|v_x(\cdot,t)\|_{L^2(\O)}=0.\ee
The combination of \eqref{x8}, \eqref{x1}, \eqref{4y}, \eqref{2.12}, and
\eqref{df8} directly yields \eqref{1y}.

 Finally, it follows from the proof in \cite{akm,ks} that  there exists some constant $c>2$
 such that  for all $(x,t)\in \ol \O\times [0,\infty)$ \bess  {c}^{-1}e^{-ct}\le
 \theta(x,t) ,\eess which together with \eqref{x5} implies that   for all $(x,t)\in
 \ol \O\times [0,\infty)$ \bess   c^{-1}e^{-cT_0}\le \theta(x,t) .\eess This
 combined with \eqref{np9}   gives \eqref{2y} provided we choose $C_4\triangleq
 \max\{C_3,c e^{ cT_0}\}.$
The proof of Lemma \ref{6y} is finished.

\end{document}